\documentclass[10pt,a4paper,reqno]{amsart} 

\usepackage{latexsym}
\usepackage{verbatim}
\usepackage{epsfig}
\usepackage{rotating}
\usepackage{amssymb}
\usepackage[T1]{fontenc}
\usepackage{afterpage}
\usepackage{color}
\usepackage{dsfont}

\usepackage{amssymb,amsfonts,amsmath,stmaryrd}

\catcode`\@=11
\def\section{\@startsection{section}{1}%
 \z@{.7\linespacing\@plus\linespacing}{.5\linespacing}%
 {\normalfont\bfseries\scshape\centering}}
\def\subsection{\@startsection{subsection}{2}%
  \z@{.5\linespacing\@plus\linespacing}{.5\linespacing}%
  {\normalfont\bfseries\scshape}}

\def\subsubsection{\@startsection{subsubsection}{3}%
 \z@{.5\linespacing\@plus\linespacing}{-.5em}
  {\normalfont\bfseries\itshape}}
\catcode`\@=12

\newtheorem{Theorem}{Theorem}
\newtheorem{Lemma}[Theorem]{Lemma}
\newtheorem{Proposition}[Theorem]{Proposition}

\def\qed{$\hfill{\vrule height 3pt width 5pt depth 2pt}$}

\newfont{\bbold}{msbm10 scaled \magstep1}
\newfont{\bbolds}{msbm7 scaled \magstep1}

\newcommand{\zs}{\mathbb{Z}}

\newcommand{\qs}{\mathbb{Q}}
\newcommand{\qss}{\mbox{\bbolds Q}}

\newcommand{\om}{\omega}

\newcommand{\GL}{\mathbb{L}}

\newcommand{\GK}{\mathbb{K}}
\newcommand{\GKb}{\overline{\mathbb{K}}}

\newcommand{\V}{\mathcal V}
\newcommand{\W}{\mathcal W}
\newcommand{\M}{\mathcal M}

\newcommand{\cS}{\mathcal S}

\newcommand{\U}{\mathcal U}
\newcommand{\C}{\mathcal C}

\newcommand{\Z}{\mathcal Z}
\newcommand{\Sn}{\mathfrak S}

\newcommand{\beq}{\begin{equation}}
\newcommand{\eeq}{\end{equation}}
\newcommand{\gf}{generating function}
\newcommand{\gfs}{generating functions}

\def\emm#1,{{\em #1}}
\newcommand{\p}{permutation}
\newcommand{\ps}{permutations}
\newcommand{\si}{\sigma}
\newcommand{\la}{\lambda}

\begin{document}
\title[Discrete excursions]
{Discrete excursions}

\author{Mireille Bousquet-M\'elou}

\address{CNRS, LaBRI, Universit\'e Bordeaux 1, 351 cours de la Lib\'eration,
  33405 Talence Cedex, France}
\email{mireille.bousquet@labri.fr}
\thanks{ } 

\begin{abstract}
It is well-known that the length \gf\ $E(t)$ of
Dyck paths (excursions with steps $\pm 1$) satisfies $1-E+t^2E^2=0$. 
The \gf\ $E^{(k)}(t)$ of
Dyck paths of height at most $k$ is $E^{(k)}=F_k/F_{k+1}$, where the $F_k$ are
polynomials in $t$ given by $F_0=F_1=1$ and $F_{k+1}=
F_k-t^2F_{k-1}$. This means that the \gf\ of these polynomials is 
$\sum_{k\ge 0} F_k z^k= 1/(1-z+t^2z^2)$.
We note that the denominator of this fraction is the minimal
polynomial of the algebraic series $E(t)$. 

This pattern persists for walks with more general steps.
For any finite set of steps $\cS$, the \gf\ $E^{(k)}(t)$ of excursions
(generalized Dyck paths) taking their steps in $\cS$ and of height at
most $k$ is the ratio $F_k/F_{k+1}$ of two polynomials. These
polynomials satisfy a linear recurrence relation with coefficients in
$\qss[t]$. Their (rational)  \gf\  can be written $\sum_{k\ge 0} F_k
z^k= N(t,z)/D(t,z)$.  The excursion \gf\ $E(t)$ is algebraic and  satisfies $D(t,E(t))=0$
(while $N(t,E(t))\not = 0$).

If $\max \cS=a$ and $\min \cS=b$, the polynomials $D(t,z)$ and
$N(t,z)$ can be taken to be respectively of degree
$d_{a,b}={{a+b}\choose a}$ and $d_{a,b} -a-b $   in $z$.
These degrees are in general optimal: for instance, when
$\cS=\{a,-b\}$ with $a $ and $b$ coprime, 
$D(t,z)$ is
irreducible, 
and is thus the minimal
polynomial of the excursion \gf\ $E(t)$.

The proofs of these results involve a slightly unusual mixture of
combinatorial and algebraic tools, among which the kernel method (which solves
certain functional  equations), 
 symmetric functions, and a pinch of Galois theory. 
\end{abstract}

\maketitle
\date{January 5th, 2007}


\section{Introduction}

\begin{figure}[bt]
     \begin{center}
      \scalebox{1}{\input{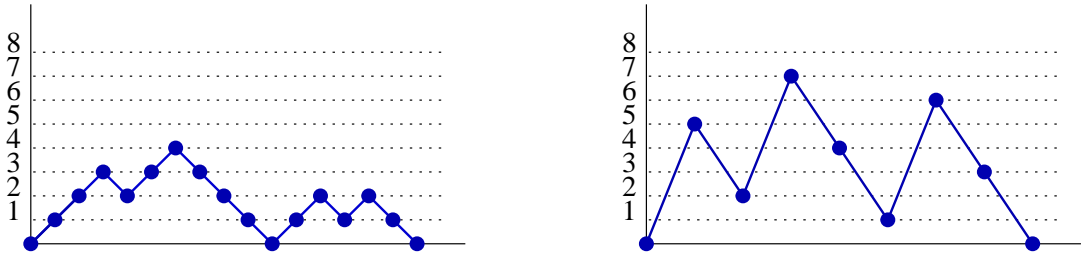}}
    \end{center}
    \caption{\emm Left:, A Dyck path of length $16$ and height 4. 
\emm Right:, An
    excursion (generalized Dyck path) of length $8$ and height 7, with
    steps in     $\cS=\{-3,5\}$.} 
\label{fig:excursion}
\end{figure}

One of the most classical combinatorial incarnation of the famous
Catalan numbers, $C_n={{2n }\choose n}/(n+1)$, is the set of \emm Dyck
paths,. These are  one-dimensional walks that start and end at $0$,
take steps $\pm 1$, and always remain at a non-negative level
(Figure~\ref{fig:excursion}, left). By factoring such walks at their first
return to 0, one easily proves that their length \gf\ $E\equiv E(t)$
is algebraic, and satisfies
$$
E=1+t^2 E^2.
$$
This immediately yields:
$$
E= \frac{1-\sqrt{1-4t^2}}{2t^2} = \sum_{n\ge 0} C_n t^{2n}.
$$
The same factorization gives a recurrence relation that defines the
 series $E^{(k)}\equiv E^{(k)}(t)$ counting Dyck paths of
height at most $k$:
$$
E^{(0)}=1 \quad \hbox{ and } 
\hbox{ for }\  k\ge1, \quad E^{(k)}=1+t ^2  E^{(k-1)} E^{(k)} .
$$
This recursion can be used to prove that $E^{(k)}$ is rational, and
more precisely, that
$$
E^{(k)}=\frac{F_k}{F_{k+1}}, \quad \hbox{ where }  \quad F_0=F_1=1
\quad \hbox{ and   } \quad F_{k+1}=F_k-t^2 F_{k-1}.
$$

The aim of this paper is to describe what happens for \emm generalized, Dyck
paths (also known as \emm excursions,) taking their steps in an
arbitrary finite set $\cS\subset \zs$ (see an example in
Figure~\ref{fig:excursion}, right). Their length \gf\ $E$ is known to be
algebraic. What is the degree of this series? How can one compute its
minimal polynomial? Furthermore, it is easy to see that the \gf \
$E^{(k)}$ can still be written $F_k/F_{k+1}$, for some polynomials
$F_k$. Does the sequence $(F_k)_k$ satisfy a linear recurrence
relation? Of what order?  How can one 
determine this  recursion? Note that any linear recursion of order $d$,
of the form
\beq\label{rec-linear}
\sum_{i=0}^d  a_iF_{k-i}=0
\eeq
with $a_i \in \qs[t]$,
gives a non-linear recursion of order $d$ for the series $E^{(k)}$,
\beq\label{non-linear}
\sum_{i=0}^d  a_i E^{(k-i+1)}\cdots E^{(k)}=0,
\eeq
and, by taking the limit $k\rightarrow \infty$, an algebraic equation
of degree $d$ satisfied by $E=\lim_k E^{(k)}$:
$$
\sum_{i=0}^d a_i E^i =0.
$$
This establishes  a close link
between the (still hypothetical) recursion for the 
sequence $F_k$ and the algebraicity of $E$. 
 The connection between~\eqref{rec-linear} and~\eqref{non-linear} is
central in the recent paper~\cite{ayyer} dealing with excursions with
steps $\pm1,\pm2$. 

\medskip

A slightly surprising outcome of this paper is that \emm
symmetric functions, are closely related to  the enumeration of excursions. 
This can be seen in the following  summary of our answers to the above
questions. 
Assume $\min\cS=-b$ and
$\max\cS=a$. Then the excursion \gf\ $E$ is algebraic of degree at
most $d_{a,b}:={{a+b}\choose a}$. The degree is exactly $d_{a,b }$ in
the \emm generic case, (to be defined), but also when $\cS=\{-b,a\}$
with $a$ and $b$ coprime. Computing a polynomial of degree $d_{a,b}$
that cancels $E$ boils down to computing the  \emm elementary
plethysms, $e_k[e_a]$ on an alphabet with $a+b$ letters, for $0\le k
\le d_{a,b}$. 

The \gf\ $E^{(k)}$ counting excursions of height at most $k$ is
rational and can be written $F_k/F_{k+1}$ for some polynomials
$F_k$. These polynomials satisfy a linear recurrence relation of the
form~\eqref{rec-linear}, of order $d_{a,b}$, which is valid for 
$k>d_{a,b}-a-b$. Moreover, $F_k$ can be expressed  as a
determinant of varying size $k$, but also  as a \emm rectangular Schur
function, taking the form of a determinant of constant size $a+b$.

\medskip
These results are detailed in the next section. Not all of them are
new.  The  \gf\ of excursions, given in
Proposition~\ref{prop:ex}, first appeared
in~\cite{bousquet-petkovsek-recurrences}, but can  be derived
from the earlier paper~\cite{gessel-factorization}. An algorithm for
computing a polynomial of degree $d_{a,b}$ that cancels $E$ was
described in~\cite{banderier-flajolet}. Hence the first part of the
next section, which deals with unbounded excursions,  is mostly a
survey (the results on the exact degree of $E$ are however new).  The
second part ---  excursions of bounded height --- 
is new, although an attempt of the same vein appears
in~\cite{banderier-these}. 

 Let us finish with the plan of this paper. 
The \emm kernel method, has become a standard tool to solve certain
 functional equations arising in various combinatorial
 problems~\cite{hexacephale,demier,prodinger}. We 
 illustrate it in Section~\ref{sec:unbounded} by counting unbounded
 excursions. We use it again in Section~\ref{sec:bounded} to obtain
 the \gf\ of excursions of bounded height. Remarkably, \emm the same result
 can be obtained by combining the transfer matrix method and the dual
 Jacobi-Trudi identity., In Section~\ref{sec:rectangle}, we determine the
 recurrence relation satisfied by the polynomials $F_k$. More
 precisely, we compute the rational series $\sum_k F_k z^k$. This is 
 equivalent to computing the  \gf\ of rectangular Schur functions
 $\sum _k s_{k^a}z^k$, where $a=\max\cS$. Finally, we discuss in
 Section~\ref{sec:degree} the exact degree of the series $E$ for
 certain step sets $\cS$. This involves a bit of Galois theory.

 \section{Statement of the results}
\label{sec:statements}
We consider  one-dimensional walks that start from $0$,
take their steps in a finite set $\cS\subset \zs$, and \emm always remain
at a non-negative level,. More formally, a (non-negative) 
walk of length $n$ will be a
sequence $(s_1, s_2, \ldots, s_n)\in \cS^n$ such that for all $i\le
n$, the partial sum $s_1+\cdots + s_i$ is non-negative. The \emm final
level, of this walk is $s_1+\cdots + s_n$, and its \emm height, is 
$\max_i s_1+\cdots + s_i$.
An \emm excursion, is a non-negative walk ending at level 0
(Figure~\ref{fig:excursion}). We are interested in the enumeration of
excursions.

The \gfs\ we consider are fairly general, in that every step $s \in \cS$ is
weighted by an element $\om_s$ of some field $\GK$ of characteristic
0. In general, we
will think of the $\om_s$ as independent indeterminates. In this case,
$\GK$ is the fraction field $\qs(\om_s, s\in \cS)$.  Of course, one
can then  specialize the $\om_s$ in various ways. 
 The length of the walks is taken
into account by an  indeterminate $t$. In particular, the
\gf\ of excursions is
$$
E:= \sum
\om_{s_1} \cdots \om_{s_n} t^n ,
$$
where the sum runs over all excursions $(s_1, s_2, \ldots, s_n)$. 
Clearly, we could incorporate the weight $t$ in the step weights
$\om_{s}$, upon replacing $\om_{s}$ by $t\om_{s}$, but we find more
convenient to keep the variable $t$  separate. 
If $\min \cS=-b$ and $\max \cS=a$, we assume
that $\om_{-b}$ and $\om_a$ are non-zero. If $d$ divides all the elements of
$\cS$, the excursion \gf\ is unchanged (up to a renaming of the weights
$\om_s$) if we  replace each $s\in \cS$ by $s/d$. Thus we can always assume that the
elements of $\cS$ are relatively 
prime.  Also, if $(s_1, s_2, \ldots, s_n)$ is an excursion, $(-s_n,
\ldots, -s_2, -s_1)$ is also an excursion, with steps in $-\cS$. Thus
the excursion series obtained for $\cS$ and $-\cS$ coincide, up to a
renaming of the weights $\om_s$.

In the expression of $E$ given below
(Proposition~\ref{prop:ex}), an important role is played by 
the following term, which encodes the steps of $\cS$:
\beq\label{P-def}
P(u)= \sum_{s\in \cS} \om_s u^s,
\eeq
where $u$ is a new indeterminate. This is a \emm Laurent, polynomial
in $u$ with coefficients in $\GK$. If $\min \cS=-b$, we define
\beq\label{K-def}
K(u)=u^b\left( 1-tP(u)\right).
\eeq
This is now a polynomial in $u$ with coefficients in 
$\GK[t]$. If $\max \cS=a$, this polynomial has degree $a+b$ in
$u$.  It has $a+b$ roots, which are fractional Laurent series (Puiseux
series) in $t$
with coefficients in $\GKb$, an algebraic closure of  
$\GK$. (We refer the reader to~\cite[Ch.~6]{stanley-vol2} for
generalities on 
the roots of a polynomial with coefficients in $\GK[t]$.) 
Exactly $b$ of these roots, say $U_1, \ldots, U_b$,
are finite at $t=0$. These roots are actually formal power series in
$t^{1/b}$. We call then the \emm small, roots of $K$.
The $a$ other roots, $U_{b+1} , \ldots ,
U_{a+b}$, are  the \emm large, roots of $K$. They are Laurent
series in $t^{1/a}$, and their first term is $c t^{-1/a}$, for some
$c\not = 0$.  
Note that $K(u)$ factors as
$$
K(u)=u^b(1-tP(u))=-t\om_{a} \prod_{i=1}^{a+b} \left(u-U_i\right),
$$
so that the elementary symmetric functions of the $U_i$'s are:
\beq\label{roots-elem}
 e_i(\U)= (-1)^i \left( \frac{\om_{a-i}}{\om_a} - \frac 1 {t\om_a} \chi_{a=i}\right),
\eeq
with $\U=(U_1, \ldots, U_{a+b})$.
We refer to~\cite[Ch.~7]{stanley-vol2} for generalities on symmetric
functions.

\subsection{Unbounded excursions}
\label{sec:intro1}
At least three different approaches have been used to count
excursions. The first one generalizes the factorization of Dyck paths
mentioned at the beginning of the introduction. It yields a system of
algebraic equations defining
$E$~\cite{Duchon98,labelle,labelle-yeh,merlini}. The 
factorization differs from one paper to another. 
To our knowledge, 
the simplest, and most systematic one, appears in~\cite{Duchon98}.

A second approach~\cite{gessel-factorization} relies on a
factorization of \emm unconstrained, walks taking their steps in
$\cS$, and on a similar factorization of formal power series. The
expression of $E$ that can be derived from~\cite{gessel-factorization}
(by combining Proposition 4.4 and the proof of Proposition 5.1)
coincides with the expression obtained by the third approach, which is
based on a step by step construction of the
walks~\cite{bousquet-petkovsek-recurrences,banderier-flajolet}. This
expression of $E$ is given in~\eqref{ex-sol} below. We repeat  in
Section~\ref{sec:unbounded} the proof of~\eqref{ex-sol} published
in~\cite{bousquet-petkovsek-recurrences},  as it 
will be extended later to count bounded excursions.
\begin{Proposition}

\label{prop:ex}
 The \gf\ of excursions is algebraic over $\GK(t)$ of degree at most
 $d_{a,b}={{a+b}\choose a}$. It can be written as:
\beq\label{ex-sol}
E= \frac{(-1)^{b+1}}{ t\om_{-b}} \prod_{i=1}^b U_i = 
\frac{(-1)^{a+1}}{ t\om_{a}} \prod_{i=b+1}^{a+b} \frac 1{U_i},
\eeq
where $U_1, \ldots, U_b$ (resp.~$U_{b+1}, \ldots, U_{a+b}$) are the
small (resp.~large) roots of the polynomial $K(u)$ given
 by~\eqref{K-def}. The quantity defined by 
\beq\label{D-expr}
D(t,z)= \prod_{I\subset  \llbracket a+b\rrbracket, \ |I|=a}
\left( 1+(-1)^a zt\om_a U_I\right),
\eeq
with
$
\llbracket a+b\rrbracket=\{1,2, \ldots, a+b\} $
 and 
$$
 U_I=\prod_{i\in I} U_i,
$$
is a polynomial in $t$ and $z$ with coefficients in $\GK$, of degree
$d_{a,b}$ in $z$, satisfying $D(t,E)=0$.
\end{Proposition}
Once the expression~\eqref{ex-sol} is established, the other
statements 
easily follow.
Indeed, the second expression of $E$ shows that
$D(t,E)=0$. Moreover, the expression of $D(t,z)$ is symmetric in the
roots $U_1, \ldots, U_{a+b}$, so that its coefficients belong to
$\GK(t)$. More precisely,  the form~\eqref{roots-elem} of
the elementary symmetric functions of the $U_i$'s shows that
$D(t,z)$ is a 
Laurent polynomial in $t$. But the valuation of $U_i$ in $t$ is at
least $-1/a$, and this implies that $D(t,z)$ is a polynomial in $t$.

Clearly, the degree of $D(t,z)$ in $z$ is $d_{a,b}={{a+b}\choose a}$. 
Thus the excursion \gf\ $E$ has degree at most $d_{a,b}$. We prove in
Section~\ref{sec:degree} that $D(t,z)$ is actually irreducible in the
two following cases:
\begin{itemize}
  \item $\cS=\llbracket -b,a\rrbracket$ and $\om_{-b},\ldots, \om_a$ are independent
  indeterminates (the \emm generic case,),
\item $\cS=\{-b,a\}$ with $\om_{-b}=\om_a=1$ and $a$ and $b$ coprime (\emm
  two-step excursions,).
\end{itemize}
As shown by Example 2 below, $D(t,z)$ is not always irreducible.

\medskip
\noindent
{\bf An algebraic equation for $E$.}
As argued above, $D(t,z)$ is a polynomial in $t$ and $z$ that
vanishes for $z=E$. However, its expression~\eqref{D-expr} involves
the series $U_i$, while  one would prefer to obtain an \emm
explicit, polynomial in $t$ and $z$.
Recall that the series $U_i$ are only known via their elementary symmetric
functions~\eqref{roots-elem}. How can one compute a polynomial
expression of $D(t,z)$? 
 The approaches based on resultants or Gröbner bases become
very quickly ineffective.

In the generic case where $\cS=\llbracket -b, a\rrbracket $ and the
weights $\om_s$ are indeterminates, $K(u)$ is the general polynomial
of degree $a+b$ in $u$, and the problem can be rephrased as
follows: Take $n=a+b$ variables $u_1, \ldots, u_n$, and expand the
polynomial
\beq\label{Q-def}
Q(z)= \prod_{|I|=a} (1-zu_I)
\eeq
in the basis of 
elementary symmetric functions of $u_1, \ldots, u_n$. For instance,
for $a=2$ and $b=1$,
\begin{eqnarray*}
Q(z)&=&(1-zu_1u_2)(1-zu_1u_3)(1-zu_2u_3)\\
&=&1-z(u_1u_2+u_1u_3+u_2u_3)+z^2(u_1^2u_2u_3+u_1u_2^2u_3+u_1u_2u_3^2)
-z^3(u_1u_2u_3)^2\\ 
&=&1-ze_2+z^2e_{3,1}-z^3e_{3,3},  
\end{eqnarray*}
while for $a=b=2$,
\beq\label{Q22}
Q(z)=1-ze_2+z^2(e_{3,1}-e_4) -z^3( e_{3,3} +e_{4,1,1}-2e_{4,2})
+z^4e_4(e_{3,1}-e_{4}) -z^5e_{4,4,2} +z^6e_{4,4,4}.
\eeq
Using the standard notation for plethysm~\cite[Appendix 2]{stanley-vol2}, 
the polynomial $Q(z)$ reads
$$
Q(z)=  \sum_{k=0} ^{d_{a,b}} (-z)^k e_k[e_a].
$$
This shows that, in the generic case, the problem of expressing
$D(t,z)$ as a polynomial in $t$ and $z$ is equivalent to expanding the
plethysms $e_k[e_a]$ in the basis of elementary symmetric functions,
for an alphabet of $n=a+b$ variables. Unfortunately, there is no
general expression for the expansion of $e_k[e_a]$ in any standard
basis of symmetric functions, and only algorithmic solutions
exist~\cite{carre,chen}.  Most of them expand 
plethysms in the basis of Schur functions. This is justified by the
representation theoretic meaning of plethysm. Still, in our
walk problem, the natural basis is that of elementary 
functions. We have used for our calculations the simple
\emm platypus,\footnote{Don't ask me why!}
algorithm presented in~\cite{banderier-flajolet}, which  only exploits
the connections between power sums and elementary symmetric
functions. This algorithm 
takes advantage automatically of  simplifications occurring in
non-generic cases. For instance, when only two steps are
allowed,  say $-b$ and $a$, all the elementary symmetric functions of
the $U_i$'s vanish, apart from $e_0(\U)$, $e_a(\U)$ and
$e_{a+b}(\U)$. It would be a shame to compute the general expansion of
$e_k[e_a]$ in the elementary basis, and \emm then, specialize most of
the $e_i$ to $0$. The platypus algorithm directly gives
the expansion of $e_k[e_a]$ modulo
the ideal generated by the $e_i$, for $ i\not = 0, a, a+b$.  For
instance, when $a=2$ and $b=1$,
$$
Q(z)\equiv 1-ze_2-z^3e_{3}^2,
$$
while for $a=2$ and $b=3$,
$$
Q(z)\equiv 1-ze_2-2z^5e_5^2+z^6e_2e_5^2-z^7e_2^2e_5^2+z^{10}e_5^4.
$$ 
and  for $a=5$ and $b=2$,
$$
Q(z)\equiv 1-ze_5 -3z^7e_7^5 +2z^8e_5e_7^5 -2z^9e_5^2e_7^5
+z^{10}e_5^3e_7^5-z^{11}e_5^4e_7^5 +3z^{14}e_7^{10}-z^{15}e_5e_7^{10}
+2z^{16}e_5^2e_7^{10}-z^{21}e_7^{15}.
$$
From the above examples, one may suspect that, in the two-step case,
the coefficient of $z^k$ in $Q(z)$ is always 
a monomial. Going back to the polynomial $D(t,z)$, and given that
$$
e_a(\U)= \frac{(-1)^{a+1}}{t\om_a} \quad \hbox{ and } \quad
e_{a+b}(\U)= (-1)^{a+b} \frac{\om_{-b}}{\om_a},
$$
this  would mean that the coefficient of $z^k$ in $D(t,z)$ is always a
monomial in $t$. This 
observation first  gave us some hope to find (in the two-step case) a
simple description of 
$D(t,z)$ and, why not, a direct combinatorial proof of
$D(t,E)=0$. However, this nice pattern does not persist:
 for $a=3$ and $b=5$, the coefficient of $z^{16}$ in $Q(z)$ contains
 $e_8^6$ and $e_3^8e_8^3$. 

\smallskip
For the sake of completeness, let us describe this platypus algorithm.
Take a polynomial
$L(z)$ of degree $n$ with constant term 1,
and define $U_1, \ldots, U_n$ implicitly by
$$
L(z)=\prod_{k=1}^n (1-zU_k).
$$
The algorithm gives a polynomial expression of
$$
Q(z)= \prod_{|I|=a} (1-zU_I) = \sum_{k=0}^d (-z)^k e_k[e_a](\U)
$$
with $d={n\choose a}$ and $\U=(U_1, \ldots, U_n)$. The only general
identity that is needed is the expansion of $e_a$ in 
power sums. This can be obtained from a series expansion  via
\beq\label{ep}
 e_a = [z^a] \exp \left(- \sum_{i\ge 1} \frac{(-z)^i}i \ p_i\right) 
= \Phi_a(p_1,\ldots, p_a) 
\eeq
for some polynomial $\Phi_a$. 
The rest of the calculation also uses series expansions, and goes as follows:
\begin{itemize}
\item
compute $p_i(\U)$ for $1\le i \le ad$ using  $p_i(\U)= i[z^i]  \log (1/L(z))$,
\item 
compute $\log  Q(z)$ up to the coefficient of $z^d$ using
\beq\label{lQexp}
  \log  Q(z)=-\sum_{i\ge 1} \frac {z^i}i\ \Phi_a(p_i(\U), p_{2i}(\U), \ldots,
  p_{ai}(\U)),
\eeq
\item compute $Q(z)$  up to the coefficient of $z^d$ using $Q(z)=\exp(\log Q(z))$.
\end{itemize}
Since $Q(z)$ has degree $d$, the calculation is complete. The
identity~\eqref{lQexp} follows from~\eqref{ep} and
$$
  \log  Q(z)=-\sum_{i\ge 1} \frac {z^i}i\  \sum_{|I|=a} U_I^i
=-\sum_{i\ge 1} \frac {z^i}i\  e_a(U_1^i, \ldots, U_n^i).
$$

  Given a set of steps $\cS$, with $\max \cS =a$, one obtains a
  polynomial expression of 
  $D(t,z)$ by applying the platypus algorithm to
$$
L(z)= \sum_{s\in \cS} \frac{\om_s}{\om_a} z^{a-s} -\frac{z^a}{t\om_a}.
$$
If the output of the algorithm is the polynomial $Q(z)$, then 
$D(t,z)=Q((-1)^{a+1}t\om_az).$

\medskip
\noindent
{\bf Example 1: Two step excursions.} The simplest walks we can
consider are obtained for  
$\cS=\{-b,a\}$ and $\om_a=\om_{-b}=1$. We always assume that $a$ and $b$
are coprime. 

If $b=1$, Proposition~\ref{prop:ex}
gives $E=U/t$, where $U$ is the only power series satisfying
$U=t(1+U^{a+1})$. Equivalently,  
$E=1+ t^{a+1} E^{a+1}$. This equation can be understood
combinatorially by looking at the first visit 
of the walk at levels $a, a-1, \ldots, 1, 0$, and factoring the walk
at these points. Of course, a similar result holds when $a=1$.

If $a, b >1$, it is still possible, but more difficult, to write
directly  a system of polynomial equations, based on  
factorizations of the walks, that define the series $E$. See for
instance~\cite{Duchon98,labelle,labelle-yeh,merlini}.  It would be
interesting to work out the precise link between the components of these
systems and the series $U_i$. To
compare both types of results, take
$a=3$ and $b=2$. On the one hand, it is shown in~\cite{Duchon98} that
$E$ is the first component of the solution of 
$$
\left\{
\begin{array}{lclllcl}
  E&=&1+ L_1R_1+L_2R_2 &  &  L_1&=&L_2R_1+L_3R_2\\
R_1&=& L_1R_2      &    &    L_2&=&L_3R_1 \\
R_2&=&tE            &  &     L_3&=&tE.
\end{array}\right.
$$
On the other hand, Proposition~\ref{prop:ex} gives $E=-U_1U_2/t$,
where $U_1, U_2$ are the small roots of $u^2=t(1+u^5)$. The platypus
algorithm gives $D(t,E)=0$ with
\beq\label{D32}
D(t,z)=1-z+t^5z^5( 2-z+z^2)+t^{10}z^{10}.
\eeq
This polynomial is irreducible. Similarly, for $a=4$ and $b=3$, $D(t,E)=0$ with
\begin{multline}\label{D43}
  D(t,z)=1-z
+ {t}^{7}z^{7} \left(5 -4\,z+z^{2}+3\,z^{3}-z^{5}+z^{6}
\right)
+ {t}^{14}z^{14} \left(
10-6\,z+3\,z^{2}+5\,z^{3}-z^{4}+z^{5} \right) \\
+{t}^{21}z^{21} \left( 10 -4\,z+3\,z^{2}+z^{3}-z^{4}\right) 
+ {t}^{28}z^{28} \left(5-z+ z^{2}-z^{3} \right)
+z^{35}{t}^{35}.
\end{multline}
We prove in Section~\ref{sec:degree} that, in the case of two step
walks,  $D(t,z)$ is always irreducible. That is, the degree of $E$ is  
exactly ${a+b}\choose a$.

\medskip\noindent
{\bf Example 2: Playing basket-ball with A. and Z.}
In a recent paper~\cite{ayyer}, the authors
 consider  excursions  with
steps in $\{\pm1, \pm2\}$, where the steps $\pm 2$ have length 2 rather
than~1. They use  factorizations of walks to  count excursions (more specifically,
excursions of bounded height).   
This problem fits in our framework by specializing the indeterminates $\om_s$
to  $\om_{-2}=\om_2=t$, and $\om_{-1}=\om_1=1$. Proposition~\ref{prop:ex} gives
$E=-U_1U_2/t^2$, 
where $U_1$ and $U_2$ are the two small roots of
$u^2=t(t+u+u^3+tu^4)$. The platypus algorithm yields
\beq\label{D-factor}
D(t,z)=\bar D(t,z) (1+t^2z)^2,
\eeq
where
\beq\label{D-basket}
 \bar D(t,z)= t^{8}{z}^{4}
-{t}^{4} \left( 1+2\,{t}^{2} \right) {z}^{3}
+{t}^{2} \left( 3+2\,{t}^{2} \right) {z}^{2}
- \left( 1+2\,{t}^{2} \right) z+1
\eeq
is the minimal polynomial of $E$. This factorization is an interesting
phenomenon, which is not
related to the unequal lengths of the steps. Indeed, the same phenomenon
occurs when $\cS=\{\pm1,\pm2\}$ and all weights are 1. In this case, one
finds:
$$
D(t,z)=\bar D(t,z) (1+tz)^2 \quad \hbox{ with } \quad
\bar D(t,z)=t^4z^4-t^2(2t+1)z^3+t(3t+2)z^2-(2t+1)z+1,
$$
so that the excursion \gf\ has degree 4 again. 

The  factorization of  $D(t,z)$ is due to the symmetry of the
set of steps. For each set of steps $\cS$ such that
$\cS=-\cS$ and weights 
$\om_s$ such that $\om_s=\om_{-s}$, the polynomial $P(u)$ given
by~\eqref{P-def} is symmetric in $u$ and $1/u$. In particular, $a=b$. 
This implies that the small and large roots of $1-tP(u)$
can be grouped by pairs: $U_{a+1}=1/U_1, \ldots, U_{2a}=1/U_a$.
In particular, if $a$ is even, the polynomial $D(t,z)$ contains
the factor $(1+t\om_az)$ at least $a\choose {a/2}$ times. In the
basket-ball case ($a=2$), this explains the factor $(1+tz)^2$
occurring in $D(t,z)$. More generally, we prove in
Section~\ref{sec:questions} that if $\cS$ is symmetric, with symmetric
weights, then the degree of $E$ is at most $2^a$, where $a=\max \cS$.

\subsection{Excursions of bounded height}
\label{sec:bounded-results}
We now turn our attention to the enumeration of excursions of height
at most $k$. These are walks  
on a finite directed graph, so that the classical transfer matrix method
applies\footnote{In language theoretic terms,  
the words of $\cS^*$ that encode these bounded excursions are
recognized by a finite automaton.}. The vertices of the graph
are $0, 1, \ldots, k$, and there is an arc from $i$ to $j$ if $j-i \in \cS$. 
The adjacency matrix of this
 graph is  $A^{(k)}=(A_{i,j})_{0\le i, j \le k}$ with
\beq\label{adjacency}
A_{i,j}=  \left\{
\begin{array}{ll}
\om_{j-i} & \hbox{if } j-i \in \cS,\\
0 & \hbox{otherwise.}  
\end{array}
\right.
\eeq
By considering the $n$th power of $A^{(k)}$, it is easy to
see~\cite[Ch.~4]{stanley-vol1} that the series $E^{(k)}$  counting
excursions of height at most $k$ is  the entry $(0,0)$
in $(1-tA^{(k)})^{-1}$. The translational invariance of our step system gives
$$
E^{(k)}= \frac{F_{k}}{F_{k+1}}
$$
where $F_0=1$ and $F_{k+1}$ is the determinant of $1-tA^{(k)}$. The
size of this matrix, $k+1$, grows with the height. 
As was already observed in~\cite{banderier-these}, the series counting 
walks confined in a strip of fixed height  can
also be expressed using determinants of size $a+b$, 
where $a=\max \cS$ and $-b=\min \cS$. However, the expressions given in the
above reference are heavy. A  different route
yields determinants that are \emm Schur functions,
in the series $U_i$ (recall that these series are the roots of the polynomial $K(u)$
given by~\eqref{K-def}).
This was shown in~\cite{mbm-ponty} for the enumeration of \emm culminating
walks,. The case of excursions is even simpler, 
as it only involves \emm rectangular, Schur functions. 

Let us recall the definition of Schur functions  in $n$ variables
$x_1, \ldots, x_n$. Let  $\delta=(n-1,n-2, \ldots, 1, 0)$.
 For any integer partition $\la$ with at most $n$ parts, $\la=(\la_1,\ldots,
 \la_n)$ with $\la_1\ge \la_2\ge \cdots\ge  \la_n\ge 0$, 
\beq\label{schur-def}
s_\la(x_1, \ldots, x_n)= \frac{a_{\delta+\la}}{a_\delta}, \quad \hbox { with  } \quad 
a_\mu= \det\left( x_i^{\mu_j}\right)_{1\le i,j\le n}.
\eeq

\begin{Proposition}\label{prop:ex-bounded}
The \gf\ of excursions of height at most $k$ is
$$
E^{(k)}= \frac{F_k}{F_{k+1}} = 
\frac{(-1)^{a+1}}{ t\om_{a}}  \frac{s_{k^a}(\U)}{s_{(k+1)^a}(\U)}
$$
where $\U=(U_1, \ldots, U_{a+b})$ is the collection of roots
of the polynomial $K(u)$ given by \eqref{K-def}, and $F_{k+1}=\det
(1-tA^{(k)})$ where $A^{(k)}$ is the adjacency matrix~\eqref{adjacency}. 
In particular,
\beq\label{F-sym}
F_k= (-1)^{k(a+1)} (t\om_a)^k s_{k^a}(\U).
\eeq
\end{Proposition}
This proposition is proved in Section~\ref{sec:bounded} in two different ways.
In Section~\ref{sec:rectangle}, we derive from the Schur expression of $F_k$ 
that  these polynomials  satisfy a linear recurrence
relation. Equivalently, the \gf\ $\sum_k F_kz^k$ is a rational function
of $t$ and $z$.
\begin{Proposition}\label{prop:link}
The \gf\ of the polynomials $F_k$ is rational, and can be written as
$$
\sum_{k\ge 0} F_k z^k= \frac{N(t,z)}{D(t,z)}
$$
where 
 $D(t,z)$ is given by~\eqref{D-expr}, and 
$N(t,z)$ has degree ${{a+b} \choose a} -a-b$ in $z$. Moreover,
$$D(t,E)=0 \quad \hbox{and} \quad N(t,E)\not = 0.$$
\end{Proposition}
In other words, the sequence $F_k$ satisfies a linear recurrence
relation of the form~\eqref{rec-linear}, of order $d_{a,b} ={{a+b}
  \choose a} $, 
valid for $k>{{a+b}   \choose a} -a-b$ (with $F_i=0$ for $i<0$). This
proposition follows from Proposition~\ref{prop:rectangle}, which 
 deals with  the \gf\ of
rectangular Schur functions of height $a$: for symmetric functions in
$n$ variables, 
\beq\label{schur-sum}
\sum_k s_{k^a} z^k=\frac{P(z)}{Q(z)}
\eeq
where $Q(z)$ is given
 by~\eqref{Q-def} and has degree $n\choose a$, while $P(z)$ has degree
 ${n\choose a} -n$.

\medskip\noindent
{\bf Computational aspects.}
We have shown in Section~\ref{sec:intro1} that, given the step set $\cS$,
the polynomial $D(t,z)$ can be computed via the platypus algorithm.
One way to determine the numerator $N(t,z)$ is to compute $F_k$
explicitly  (e.g. as the determinant of $(1-tA^{(k)})$) for $k\le
\delta:={{a+b}   \choose a} -a-b$, and then to compute $N(t,z)=D(t,z)\sum_k F_k
z^k$  up to the coefficient of $z^\delta$.
 
In the generic case, computing the \gf\ of the polynomials $F_k$
 boils down to computing the \gf\ 
~\eqref{schur-sum}. Again,
 the platypus algorithm can be used to 
 determine $Q(z)$ in terms of the elementary symmetric functions. In
 order to determine $P(z)$, we express the Schur functions
 $ s_{k^a}$, for $k\le \delta:=
{{n}   \choose a} -n$, in the elementary basis. This can be done using
the dual Jacobi-Trudi identity (see Section~\ref{sec:rectangle} for
details). 
One finally obtains $P(z)$  by expanding the product $Q(z) \sum_k
s_{k^a} z^k$ in the elementary 
 basis up to order $\delta$. 
For instance, for $a=b=2$, 
$$
\sum_{k\ge 0} s_{k^a} z^k= \frac{1-z^2e_4}{Q(z)},
$$
where $Q(z)$ is given by~\eqref{Q22}. More values of $P(z)$ are given
in Section~\ref{sec:P}. Let us now revisit the examples
of Section~\ref{sec:intro1}.
 
\medskip\noindent
{\bf Example 1: Two step excursions.} When $\cS=\{a,-1\}$,
one has $D(t,z)=1-z+t^{a+1}z^{a+1}$.  
The polynomials $F_k$ satisfy the recursion
$F_k=F_{k-1}-t^{a+1}F_{k-a-1}$, which can be 
understood combinatorially using Viennot's theory of \emm heaps of
pieces,~\cite{viennot-heaps}. 
Via this theory, $F_k$ appears as the \gf\ of trivial heaps of
segments of length $a$ on the 
line $\llbracket 0, k \rrbracket$, each segment being weighted by
$-t^{a+1}$. The recursion is valid for $k\ge 1$, with $F_0=1$ and
$F_i=0$ for $i<0$. The \gf\ of the $F_k$'s is
$$
\sum_{k\ge 0} F_k z^k =\frac 1 {1-z+t^{a+1}z^{a+1}}.
$$
When $a=3$ and $b=2$, the minimal polynomial of the excursion series
$E$ is given by~\eqref{D32} 
and the \gf\ of the polynomials $F_k$ is found to be
$$
\sum_{k\ge 0} F_k z^k =\frac {1+t^5z^5} {1-z+t^5z^5( 2-z+z^2)+t^{10}z^{10}}.
$$
For $a=4$ and $b=3$, we refer to~\eqref{D43} for the minimal polynomial
$D(t,z)$ of $E$, and 
$$
\sum_{k\ge 0} F_k z^k =\frac {
1+{t}^{7}{z}^{7} \left( 4+{z}^{3}+{z}^{4} \right) 
+{t}^{14}{z}^{14} \left( 6+{z}^{3} \right) 
+4\,{t}^{21}{z}^{21}
+{t}^{28}{z}^{28}
}
 {D(t,z)}.
$$

\medskip\noindent
{\bf Example 2: Basket-ball again.} 
For $\cS=\{\pm1,\pm2\}$ with $\om_{-2}=\om_2=t, \om_{-1}=\om_1=1$,  
$$
\sum_{k\ge 0} F_k z^k =
{\frac {1-{t}^{2}z}{ \left( 1+{t}^{2}z \right)  \left( 
1-z (1+ 2{t}^{2})+\,{z}^{2}{t}^{2}(3+2{t}^{2})
-{z}^{3}{t}^{4}(1+2\,{t}^{2} )
+{z}^{4}{t}^{8}
\right) }}.
%
$$
The denominator is not irreducible. Its second factor is the minimal
polynomial of $E$, 
see~\eqref{D-basket}.  Moreover, comparing to~\eqref{D-factor} shows
that  $N(t,z)$ and $D(t,z)$ have a factor $(1+t^2z)$ in common. 
A similar phenomenon occurs for
$\cS=\{\pm1,\pm2\}$ but now $\om_s=1$ for all $s$. In this case,
$$
\sum_{k\ge 0} F_k z^k =
{\frac {1-tz}{ \left(1+ zt \right)  \left( 
1-z \left( 1+2\,t \right) 
+ t\left( 2+3\,{t} \right) {z}^{2}
- {t}^{2} \left(1+ 2\,{t}\right) {z}^{3}
+{z}^{4}{t}^{4}
\right) }}.
$$
Again, the minimal polynomial of $E$ is the second factor of the
denominator, and $N(t,z)$ and $D(t,z)$ have a factor $(1+tz)$ in
common.

\section{Enumeration of unbounded excursions}
\label{sec:unbounded}
Here we establish the expression~\eqref{ex-sol} of the excursion \gf\
$E$. The proof is based on a step-by step construction of
non-negative walks with steps in $\cS$, and on the so-called \emm
kernel method,. This type of argument is by no means original. The 
proof that we are going to present can be found
in~\cite[Example~3]{bousquet-petkovsek-recurrences}, then
in~\cite{banderier-flajolet}, and finds its origin
in~\cite[Ex.~2.2.1.4 and 2.2.1.11]{knuth}. 
The reason why we repeat the proof is  because it
will be adapted  in Section~\ref{sec:bounded} to count excursions of
bounded height.

Let $\W$ be the set of walks that start from 0, take their steps in
$\cS$, and always remain at a non-negative level. Let $W(t,u)$
be their \gf, where the variable $t$ counts the length,  the
variable $u$ counts the final height, and each step $s \in \cS$ is
weighted by  $\om_s$:
$$
W(t,u)= \sum _{(s_1, s_2, \ldots, s_n)\in \W} \om_{s_1} \cdots \om_{s_n}
t^n u^{s_1+ \cdots + s_n}.
$$
We often denote
$W(t,u)\equiv W(u)$, and use the notation $W_h$ for the \gf\ of walks
of $\W$ ending at height $h$:
$$
W(t,u)=\sum_{h\ge 0} u^h W_h \quad \hbox{where} \quad
W_h=\sum _{(s_1, s_2, \ldots, s_n)\in \W \atop
s_1+  \cdots + s_n =h} \om_{s_1} \cdots \om_{s_n}
t^n.
$$
 A non-empty walk of $\W$ is obtained by adding a step of $\cS$ at
the end of another walk of $\W$. However,  we must avoid
adding a step $i$ to a walk ending at height $j$, if
$i+j<0$. This gives
$$
W(u)=1+t \left(\sum_{s\in \cS} \om_s u^s\right) W(u) 
-t \sum_{i \in \cS, j\ge 0\atop i+j<0} \om_i u^{i+j} W_j.
$$
Let $\min \cS =-b$.  Rewrite the above equation so as to involve  only
non-negative powers of $u$:
\beq\label{ex-main-eq}
u^b(1-tP(u)) W(u)=u^b-t\sum_{h=1}^b u^{b-h} \sum_{i \in \cS, j\ge 0\atop i+j=-h }
 \om_i W_j,
\eeq
with
$
P(u)= \sum_{s\in \cS} \om_s u^s.
$
The coefficient of $W(u)$ is the \emm kernel,  $K(u)$ of the equation,
given in~\eqref{K-def}.
As above, we denote by $U_1, \ldots, U_b$ (resp. $U_{b+1}
, \ldots , U_{a+b}$) the roots of $K(u)$ that are
 finite (resp. infinite) at $t=0$.
For $1\le i \le b$, the series $W(U_i)$ is well-defined (it is a
formal power series in $t^{1/b}$). 
The left-hand side of~\eqref{ex-main-eq} vanishes for
$u=U_i$, with $i\le b$, and so the right-hand side vanishes too. But
the right-hand 
side is a polynomial in $u$, of degree $b$, leading coefficient 1,
and it vanishes at $u=U_1, \ldots, U_b$. This gives
$$
u^b(1-tP(u)) W(u)= \prod_{i=1}^b (u-U_i).
$$
As the coefficient of $u^0$ in the kernel is $-t\om_{-b}$, setting $u=0$
in the above equation gives the \gf\ of excursions:
$$
E=W(0)= \frac{(-1)^{b+1}}{t\om_{-b}} \prod_{i=1}^b U_i.
$$
This is the first expression in~\eqref{ex-sol}. The second follows
using 
\beq\label{product-roots}
U_1 \cdots U_{a+b}= (-1)^{a+b} \om_{-b}/\om_a
\eeq
(see~\eqref{roots-elem}).
\qed

\medskip 
\noindent
{\bf Remark.}
There exists an alternative way to solve~\eqref{ex-main-eq}, which does not
exploit the fact that the  
right-hand side of~\eqref{ex-main-eq} has degree $b$ in $u$. This
variant will be useful in the enumeration of bounded excursions. Write
$$
Z_{-h}=\sum_{i \in \cS, j\ge 0\atop i+j=-h } \om_i W_j,
$$
so that  the right-hand side of~\eqref{ex-main-eq} reads
$$
u^b -t\sum_{h=1}^b u^{b-h} Z_{-h}.
$$
This term vanishes for $u=U_1,
\ldots, U_b$. Hence the $b$ series $Z_{-1}, \ldots, Z_{-b}$ satisfy
the following system of $b$ linear equations: For $U=U_i$, with
$1\le i \le b$,
$$
 \sum_{h=1}^{b} U^{b-h} Z_{-h} = U^{b}/t.
$$
In matrix form, we have $\M\Z=\C/t$, where $\M$ is the square matrix
of size $b$ given by
$$
\M=\left(
\begin{array}{cccccccccccccccccc}
  U_1^{b-1}  & U_1^{b-2} & \cdots &  U_1^{1}  & 1 \\
  U_2^{b-1}  & U_2^{b-2} & \cdots &  U_2^{1}  & 1 \\
\vdots &&&&\vdots \\
 U_b^{b-1}  & U_b^{b-2} & \cdots &  U_b^{1}  & 1
\end{array}
\right),
$$
$\Z$ is the column vector  $(Z_{-1} , \ldots, Z_{-b})$,
and $\C$ is the column vector $(U_1^{b},\ldots, U_{b}^{b})$.
The determinant of $\M$ is the Vandermonde in $U_1, \ldots, U_b$, and
it is non-zero because the $U_i$ are distinct. We are especially
interested in the unknown $Z_{-b}=\om_{-b} E$. Applying Cramer's rule to
solve the above system yields
$$
Z_{-b}= \frac{(-1)^{b+1}}  t  \frac 
{\det (U_i^{b-j+1})_{1\le   i,j\le b}} 
{\det (U_i^{b-j})_{1\le   i,j\le b}} .
$$
The two determinants coincide, up to a factor $U_1\ldots U_b$, and we
finally obtain
$$
E= \frac {Z_{-b}}{\om_{-b}}= \frac{(-1)^{b+1}}  {t\om_{-b}} U_1\cdots U_b.
$$
\section{Enumeration of bounded excursions}
\label{sec:bounded}
As argued in Section~\ref{sec:bounded-results}, the \gf\ of excursions
of height at most $k$ is
\beq\label{Ek-F}
E^{(k)}=\frac {F_k}{F_{k+1}},
\eeq
where $F_{k+1}=\det(1-tA^{(k)})$ and $A^{(k)}$ is the adjacency
matrix~\eqref{adjacency} 
describing the allowed  steps in the
interval $\llbracket 0,k\rrbracket$. In order to prove
Proposition~\ref{prop:ex-bounded}, it remains to establish the
expression~\eqref{F-sym} of the polynomial $F_k$ as a Schur
function of $U_1, \ldots, U_{a+b}$. We give two proofs. The first one uses the
dual Jacobi-Trudi identity to identify $F_k$ as a Schur function. The
second determines $E^{(k)}$ in terms of Schur functions via the kernel
method, and the Schur expression of $F_k$ then follows from~
\eqref{Ek-F} by induction on $k$ (given that $F_0=1$).

\medskip
\noindent
{\bf First proof via the Jacobi-Trudi identity.}
The dual Jacobi-Trudi identity expresses
Schur functions as a determinant in the elementary symmetric
functions $e_i$~\cite[Cor.~7.16.2]{stanley-vol2}:  for any partition
$\la$, 
$$
s_\la= \det\left( e_{\la'_j+i-j}\right)_{1\le i,j\le \la_1},
$$
where $\la'$ is the conjugate of $\la$. 
Apply this identity  to $\la=(k+1)^a$. Then $\la'=a^{k+1}$ and
$$
s_{(k+1)^a}=\det J^{(k)} \quad \hbox{ with }
\quad J^{(k)}=( e_{a+i-j})_{1\le i, j\le k+1}.
$$
Now, specialize this to symmetric functions in  the $a+b$ variables
$\V=(V_1, \ldots, V_{a+b})$ where  $V_i=-U_i$ for all $i$.
  By~\eqref{roots-elem},  the elementary symmetric functions of the $V_i$
are 
$$
e_i(\V)=\frac{\om_{a-i}}{\om_a} - \frac 1 {t\om_a} \chi_{i=a}= -\frac
1{t\om_a}
\left(  \chi_{i=a}-t \om_{a-i}\right).
$$ 
This shows that the matrix $J^{(k)}$ coincides with
$-(1-tA^{(k)})/(t\om_a)$, so that  
%
$$
s_\la(\V)= (-t\om_a)^{-(k+1)} F_{k+1}=(-1)^{a(k+1)}s_\la(\U),
$$
since $s_\la$ is homogeneous of degree $a(k+1)$. This gives the Schur
expression of $F_{k+1}$. 
\qed

\medskip
\noindent
{\bf Second proof via the kernel method.} We adapt the step by step approach of
Section~\ref{sec:unbounded} 
to count excursions of height at most $k$. Let $W^{(k)}(t,u)\equiv
W^{(k)}(u)$ be the \gf\ of non-negative walks of height at most $k$.
 As before, we count them by their length (variable $t$) and final height
($u$) with  multiplicative weights $\om_s$ on the steps. We
use notations similar to those of Section~\ref{sec:unbounded}. When
constructing walks step by step, we must still avoid going below level
0, but also above level $k$. This yields:
$$
W^{(k)}(u)=1+t \left(\sum_{s\in \cS} \om_s u^s\right) W^{(k)}(u) -
t \sum_{i \in \cS, j\ge 0 \atop  i+j>k \hbox{\tiny{ or} } i+j <0} \om_i u^{i+j} W^{(k)}_j,
$$
or, with $\min \cS =-b$,
\beq\label{ex-bounded-main-eq}
u^b(1-tP(u)) W^{(k)}(u)=u^b-t\sum_{h=k+1}^{k+a} u^{b+h} Z^{(k)}_h
- t\sum_{h=1}^b u^{b-h} Z^{(k)}_{-h},
\eeq
where
$$
Z_h^{(k)}= \sum_{i \in \cS, j\ge 0\atop i+j=h }  \om_i W^{(k)}_j.
$$
The series $W^{(k)}(u)$ is now a \emm polynomial, in $u$ (with
coefficients in the ring of power series in $t$). This implies that
any root $U_i$ of the kernel 
$K(u)=u^b(1-tP(u))$ can be legally substituted for $u$
in~\eqref{ex-bounded-main-eq}.  The right-hand side and the left-hand
side then vanish, and provide a system of $a+b$ linear equations
satisfied by the $Z_h$:
 For $U=U_i$, with $1\le i \le a+b$,
$$
 \sum_{h=k+1}^{k+a} U ^{b+h} Z^{(k)}_h+ \sum_{h=1}^{b} U^{b-h} Z^{(k)}_{-h} = U^{b}/t.
$$
In matrix form, we have $\M^{(k)}\Z^{(k)}=\C/t$, where $\M^{(k)}$ is
the square matrix of size $a+b$ given by
$$
\M^{(k)}=\left(
\begin{array}{cccccccccccccccccc}
 U_1^{a+b+k}  & U_1^{a+b+k-1} & \cdots & U_1^{b+k+1} & 
U_1^{b-1} & U_1^{b-2} & \cdots & 1 \\
 U_2^{a+b+k}  &\cdots &&&&& \cdots & 1 \\
\vdots &&&&&&&\vdots \\
 U_{a+b}^{a+b+k}  & U_{a+b}^{a+b+k-1} & \cdots & U_{a+b}^{b+k+1} & 
U_{a+b}^{b-1} & U_{a+b}^{b-2} & \cdots & 1 \\
\end{array}
\right),
$$
$\Z^{(k)}$ is the column vector  $(Z^{(k)}_{k+a} ,\ldots,
Z^{(k)}_{k+1}, Z^{(k)}_{-1} , \ldots, Z^{(k)}_{-b})$, 
and $\C$ is the column vector $(U_1^{b},\ldots, U_{a+b}^{b})$.
 We are especially interested in the series $Z^{(k)}_{-b}=\om_{-b} E^{(k)}$. 
Cramer's rule now gives
\beq\label{Z-b}
Z^{(k)}_{-b}= \frac {(-1)^{b+1}} t \frac {\det(U_i^{a+b+k}, \ldots,
  U_i^{b+k+1}, U_i^{b},   U_i^{b-1} ,
\ldots, U_i )_{1\le i \le a+b}} {\det \M^{(k)}},
\eeq
provided $\det \M^{(k)}\not = 0$.
In view of the definition~\eqref{schur-def} of Schur functions, this yields:
$$
E^{(k)}= \frac{Z^{(k)}_{-b}}{\om_{-b}}= \frac {(-1)^{b+1}}
{t\om_{-b}}\ 
U_1 \cdots U_{a+b} \  \frac{s_{k^a}(\U)}{s_{(k+1)^a}(\U)} .
$$
Thanks to~\eqref{product-roots},  the 
\gf\ of excursions of height at most $k$ can finally be rewritten 
$$
E^{(k)}= 
\frac{(-1)^{a+1}}{ t\om_{a}}  \frac{s_{k^a}(\U)}{s_{(k+1)^a}(\U)}.
$$
Using~\eqref{Ek-F}, we finally express the polynomial $F_k$ in terms
of Schur functions: 
$$
F_{k} = \frac 1 {E^{(0)} \cdots E^{(k-1)} }  
=(-1)^{k(a+1)} (t\om_a)^k s_{k^a}.
$$
We still have  to prove that
the determinant of $\M^{(k)}$ is non-zero. Whether $ \M^{(k)}$ is
singular or not,  the following variant of~\eqref{Z-b} remains valid:
\begin{eqnarray*}
  {\det \M^{(k)}} Z^{(k)}_{-b}&= &
\frac {(-1)^{b+1}} t \ {\det(U_i^{a+b+k}, \ldots,
  U_i^{b+k+1}, U_i^{b},   U_i^{b-1} , \ldots, U_i )_{1\le i \le a+b}}  \\
&= &\frac {(-1)^{b+1}} t \ V(\U)\ 
s_{k^a}(\U) \ \prod_{i=1}^{a+b} U_i,
\end{eqnarray*}
where $V(\U)$ denotes the Vandermonde in the $U_i$'s.  Since these
series are distinct and non-zero, this shows  
that if $\det\M^{(k)}=0$, that is,  $s_{(k+1)^a}(\U)=0$, then  $s_{k^a}(\U)=0$ as
    well. But this would finally imply $s_0(\U)=0$, while $s_0(\U)=1$.
Thus  $\det\M^{(k)}\not = 0$, and the second proof of
Proposition~\ref{prop:ex-bounded} is now complete.
\qed


\section{Generating functions of rectangular Schur functions}
\label{sec:rectangle}
We will now prove Proposition~\ref{prop:link}, which connects the
(algebraic) excursion \gf\ $E$ to the polynomials $F_k$ occurring in
the (rational) \gf\ $E^{(k)}$ counting excursions of height at most
$k$. Now that we have expressed $F_k$ as a Schur
function~\eqref{F-sym}, Proposition~\ref{prop:link} is a specialization of
the following result. 

\begin{Proposition}\label{prop:rectangle}
Let $1\le a\le n$. The \gf\ of
rectangular Schur functions of length $a$  in $n$ variables 
$u_1, \ldots, u_n$ is
$$
\sum_{k\ge 0}s_{k^a} z^k = \frac{P(z)}{Q(z)}
$$
where
\beq\label{Q-def1}
Q(z)=  \prod_{I\subset  \llbracket n\rrbracket, \ |I|=a}
\left( 1- z u_I\right)= \sum_{k\ge0} (-1)^k z^ke_k[e_a]
\eeq
has degree $n\choose a$ in $z$ and $P(z)$ has degree ${n\choose a}-n$.
(We have used the notation
$
u_I= \prod_{i\in I} u_i.
$)
Moreover, for all $J$ of cardinality $a$, 
\beq\label{PuJ}
P(1/u_J)= \prod_{I : |I|=a, |I\Delta J|\ge 4} (1-u_I/u_J).
\eeq
\end{Proposition}

\begin{proof}
Let us write $n=a+b$.  By definition of Schur functions,
\beq\label{schur-rectangle}
s_{k^a}= \frac 1 {V_n} \det\left( ( u_i^{n+k-1}, \cdots, u_i^{b+k},
u_i^{b-1}, \cdots, 1)_{1\le i \le n}\right),
\eeq
where $V_n=\prod_{1\le i<j\le n} (u_i-u_j)$. Thus
\begin{eqnarray}
  \sum_{k\ge 0}s_{k^a} z^k &= &
\frac 1 {V_n} \sum_{k\ge 0} z^k \sum_{\si \in \Sn_{n}} \varepsilon(\si) \ 
\si \left( 
u_1^{n+k-1}\cdots u_a^{b+k}u_{a+1}^{b-1}\cdots 
u_{n-1}^{1}u_n^0
\right) \label{schur-series1}\\
&=& \frac 1 {V_n} \sum_{\si \in \Sn_n} \varepsilon(\si) \ \si\left( 
\frac{u_1 ^{n-1} \cdots u_a^b u_{a+1}^{b-1} \cdots u_{n-1}^1 u_n^0}
{1-zu_1\cdots u_a}
\right),\nonumber
\end{eqnarray}
where $\si$ acts on functions of $u_1, \ldots, u_n$ by permuting the
variables:
$$
\si F(u_1, \ldots, u_n)= F(u_{\si(1)}, \ldots, u_{\si(n)}).
$$
Equivalently,
$$
 \sum_{k\ge 0}s_{k^a} z^k = \frac{P(z)} {Q(z) }
$$
where $Q(z)$ is given by~\eqref{Q-def1} and
\beq\label{P-expr}
P(z)=  \frac 1 {V_n} \sum_{\si \in \Sn_n} \varepsilon(\si) \ \si\left( 
u_1 ^{n-1} \cdots u_n^0 
 \prod_{
|I|=a, I \not = \llbracket  a\rrbracket}
\left( 1- z u_I
\right)\right).
\eeq
The above expression suggests that the degree of $P(z)$ could be as large as
${n\choose a } -1$, while we claim it is only ${n\choose a}-n$.
To explain this gap, it suffices to notice that the
determinant~\eqref{schur-rectangle} 
vanishes for $k \in \{-n+1, -n+2, \ldots, -1\}$. Thus the sum over $k$
in~\eqref{schur-series1} could just as well start at $k=-n+1$, giving:
 $$
  z^{n-1}\sum_{k\ge 0}s_{k^a} z^k =
 \frac 1 {V_n} \sum_{\si \in \Sn_n} \varepsilon(\si) \ \si\left( 
\frac{u_1 ^{0}u_2^{-1} \cdots u_a^{-a+1} u_{a+1}^{b-1} \cdots u_{n-1}^1 u_n^0}
{1-zu_1\cdots u_a}
\right).
$$
This provides the following alternative expression of $P(z)$:
\beq\label{P-expr1}
z^{n-1} P(z)=  \frac 1 {V_n} \sum_{\si \in \Sn_n} \varepsilon(\si) \ \si\left( 
u_1 ^{0}u_2^{-1} \cdots u_a^{-a+1} u_{a+1}^{b-1} \cdots u_{n-1}^1 u_n^0
 \prod_{
|I|=a, I \not = \llbracket
 a\rrbracket}
\left( 1- z u_I
\right).
\right)
\eeq
The right-hand side is a polynomial in $z$ of degree (at most)
${n\choose a } -1$, and this polynomial  is the product of $P(z)$ and
$z^{n-1}$. This 
shows that $P(z)$ has degree at most ${n\choose a}-n$. Moreover, by
extracting the coefficient of $z^{{n\choose a } -1}$ in the above
identity, one finds:
$$
[z^{{n\choose a}-n}] P(z) = \frac 1 {V_n} \sum_{\si \in \Sn_n} \varepsilon(\si) \ \si\left( 
u_1 ^{0}u_2^{-1} \cdots u_a^{-a+1} u_{a+1}^{b-1} \cdots u_{n-1}^1 u_n^0
 \prod_{
|I|=a, I \not = \llbracket
 a\rrbracket}
\left( -  u_I \right).\right)
$$
Up to a sign and a power of $u_1\cdots u_n$, the sum over $\sigma$ is
the Vandermonde  in the $u_i$'s. Finally,
\beq\label{P-dominant}
[z^{{n\choose a}-n}] P(z)= (-1)^{{n\choose a}+ab-1} \left(u_1\cdots
u_n\right)
^{{{n-1}\choose {a-1}}-a},
\eeq
so that $P(z)$ has degree ${n\choose a}-n$ exactly.

It remains to determine $P(1/u_J)$, for $|J|=a$. We specialize the
expression~\eqref{P-expr} of $P(z)$ to the case $z=1/u_J$. The only \ps\
$\sigma$ having a non-zero contribution are those such that
$\sigma(\llbracket a\rrbracket )= J$. Every such \p\ $\sigma$ can be
written in a unique way $\sigma=\pi \tau\sigma_J$, where $\sigma_J$ is
the shortest \p\ sending $\llbracket a\rrbracket$ to $J$, and $\tau$
(resp.~$\pi$) is any \p\ on $J$ (resp.~$^cJ$). Thus, if $J=\{j_1,
\ldots, j_a\}$ with $j_1<\ldots< j_a$ and $^cJ=\{k_1, \ldots, k_b\}$
with $k_1< \ldots< k_b$, we have
\begin{eqnarray*}
  P(1/u_J)&=& \prod_{|I|=a, I\not = J} \left(1-u_I/u_J\right)
\frac {\varepsilon(\si_J) } {V_n} \sum_{\tau  \in \Sn(J)}  \varepsilon(\tau) \tau
\left( u_{j_1}^{n-1} \cdots u_{j_a}^b\right) 
\sum_{\pi \in \Sn(^cJ)}  \varepsilon(\pi) \pi
\left( u_{k_1}^{b-1} \cdots u_{k_b}^0\right) 
\\
&=&\prod_{|I|=a, I\not = J} \left(1-u_I/u_J\right)\frac
  {\varepsilon(\si_J) } {V_n} \  u_J^b \ V(J)\  V(^cJ),
\end{eqnarray*}
where $V(J)$ denotes the Vandermonde  in the variables
$u_j, j\in J$. This is easily seen to be equivalent to~\eqref{PuJ}.

\end{proof}

We can now complete the proof of Proposition~\ref{prop:link}. 
We combine the Schur expression of
$F_k$ given in Proposition~\ref{prop:ex-bounded} with
Proposition~\ref{prop:rectangle}. Set $n=a+b$. The indeterminates $u_1, \ldots,
u_{a+b}$ are specialized to $U_1, \ldots, U_{a+b}$, and we obtain:
$$
\sum_{k\ge 0}F_kz^k= \frac{P((-1)^{a+1}t\om_a z)}{Q((-1)^{a+1}t\om_a
  z)}= \frac{N(t,z)}{D(t,z)}
$$
where $D(t,z)=Q((-1)^{a+1}t\om_a  z)$ is exactly the
polynomial~\eqref{D-expr}. The dominant coefficient of $P(z)$, given
by~\eqref{P-dominant}, does not vanish when specializing
$u_i$ to $U_i$. Thus $N(t,z)=P((-1)^{a+1}t\om_a  z)$  has degree ${n\choose
  a}-n$ exactly. We have already seen that the excursion \gf\ $E$ given in
Proposition~\ref{prop:ex} satisfies $D(t,E)=0$. Now, since
$E=(-1)^{a+1}/(t\om_a U_J)$, with $J=\{b+1, \ldots, a+b\}$, 
$$
N(t,E)= P(1/U_J)= \prod_{|I|=a, |I\Delta J|\ge 4} \left(1-U_I/U_J\right)
$$
by~\eqref{PuJ}. Recall that $U_{b+1}, \ldots, U_{a+b}$ are  the
roots of $K(u)$ with valuation $-1/a$, while the $b$ other roots have
valuation $1/b$. This implies that $U_I\not = U_J$ for $I\not = J$, so
that $N(t,E)\not = 0$.
\qed


\section{The degree of the excursion \gf}
\label{sec:degree}
We conclude this paper by proving that the results stated in
Section~\ref{sec:statements} 
are, in a sense, optimal. We have defined in~\eqref{D-expr} a polynomial
$D(t,z)$, of degree $d_{a,b}= {{a+b}\choose a}$, which satisfies
$D(t,E)=0$ and is the  denominator of the 
rational series $\sum_k F_k z^k$. We prove that $D(t,z)$ is
irreducible in the following two cases: 
\begin{itemize}
  \item $\cS=\llbracket -b,a\rrbracket$ and $\om_{-b},\ldots, \om_a$ are independent
  indeterminates,
\item $\cS=\{-b,a\}$ with $\om_{-b}=\om_a=1$ and $a$ and $b$ coprime.
\end{itemize}
In the first case, the kernel $K(u)$ is essentially the general
algebraic equation of degree $a+b$, so that the result may be
predictable. The idea is that there are no non-trivial relations
between the series 
$U_i$. The second case is less obvious. 
\begin{Proposition}\label{prop:degree}
  In the above two cases, the \gf\ of excursions with steps in $\cS$
is algebraic of degree $d_{a,b}= {{a+b} \choose
  a}$. Its minimal polynomial is given by~\eqref{D-expr}.
\end{Proposition}
Recall, from Example~2 in
Section~\ref{sec:intro1}, that $E$ has sometimes degree less than
$d_{a,b}$ (for instance when $\cS=\{\pm1, \pm2\}$ with weights 1). 

 The key tool is the study of the Galois group of the polynomial $K(u)$.
We begin with a  condition implying the irreducibility  of $D(t,z)$.
\begin{Lemma}
   Let $\cS$ be a finite set of steps with weights $\om_s \in \GK$. Let
  $a=\max \cS$, $-b=\min \cS$ and $n=a+b$. Let $K(u)$ be the
  polynomial in $u$, with coefficients in $\GK(t)$, defined by
  \eqref{K-def}.  

If the Galois group of $K(u)$ over $\GK(t)$, seen as
  a permutation group of the $U_i's$, is the
  full symmetric group $\Sn_n$, then the product $U_1 \ldots U_b$ of the
small roots of $K(u)$ has degree  $d_{a,b}= {{a+b} \choose
  a}$. In other words, the polynomial $D(t,z)$ given by~\eqref{D-expr}
  is irreducible.
\end{Lemma}
\begin{proof}
The extension $\GK(t,U_1, \ldots,
U_n)$ of $\GK(t)$ is normal by construction, and separable since we
have assumed $\GK$ to be of characteristic 0.  Assume that the Galois
group of  $\GK(t,U_1, \ldots,
U_n)$ over $\GK(t)$ is $\Sn_n$.
By the main result of Galois theory, the correspondence  $\Phi$ between
subgroups $G$ of $\Sn_n$ and sub-extensions $\GL$ of $\GK(t,U_1, \ldots,
U_n)$  defined by
$$
\Phi(G)
=\GL
= \{ x \in \GK(t,U_1, \ldots, U_n): \ \si(x)=x \hbox{ for all } \si \in
G\}
$$
is bijective. Its inverse is given by
$$
\Phi^{-1}(\GL)=G
= \{\si \in \Sn_n: \si(x)=x \hbox{ for all }  x \in \GL\}.
$$
Moreover, the degree of $\GK(t,U_1, \ldots, U_n)$ over $\GL$ is $|G|$.

In particular, let $\GL=\GK(t,U_1\cdots U_b)$ be the extension of
$\GK(t)$ generated by the product of the small roots. Given that $U_1,
\ldots, U_b$ have valuation $1/b$ in $t$, while $U_{b+1}, \ldots,
U_{a+b}$ have valuation $-1/a$, the only permutations $\si$ of $\Sn_n$
that leave $U_1\ldots, U_b$ unchanged are those that fix the set
$\llbracket b \rrbracket$. That is,  $\Phi^{-1}(\GL)\simeq\Sn_b\times
\Sn_a$. Thus $\GK(t, U_1, \ldots, U_n)$ has degree $a!b!$ over $\GL$,
degree $(a+b)!$ over $\GK(t)$, so that $\GL=\GK(t,U_1\cdots U_b)$ has
degree $ {{a+b} \choose  a}$ over $\GK(t)$.
\end{proof}

We now apply the above lemma to prove Proposition~\ref{prop:degree}.

\noindent{\em Proof of Proposition~\ref{prop:degree}.}
In the first case, $K(u)$ is  the general equation of degree $n=a+b$. It
is well-known that its Galois group is $\Sn_n$. See for
instance~\cite{stewart}. 

In the second case, we want to prove that the Galois group of
$K(u)=u^b-t(1+u^{a+b})$ over $\qs(t)$ is $\Sn_n$, with $n=a+b$. 
This has been proved for trinomials $u^{a+b}+\alpha u^b+ \beta$
with \emm two indeterminate, coefficients $\alpha$
and $\beta$ (see~\cite{smith,cohen}), and for  some trinomials
with \emm rational, coefficients~\cite{osada,salinier}. The latter results
are of course harder that the former. Given that we could not find any reference
dealing with trinomials involving exactly \emm one, indeterminate
coefficient, we will  rely on the strong results obtained for
trinomials of $\qs[u]$.

We first note that it suffices to prove that the trinomial
$u^b-t_0(1+u^{a+b})$ has Galois group $\Sn_n$ over $\qs$ for some
rational number $t_0$. Since $a$ and $b$ are coprime, Theorem~8
of~\cite{schinzel} implies that there exist only finitely many $\alpha\in
\zs$ such that $u^{a+b}+\alpha u^b+1$ is reducible. Thus we can choose
 $\alpha\in \zs$, coprime with $n=a+b$, and such that the above
trinomial is irreducible. Then by~\cite[Thm.~1]{osada}, this trinomial has
Galois group $\Sn_n$ over $\qs$.
\qed

\section{Concluding remarks and questions}
\label{sec:questions}

\subsection{The degree of the excursion \gf}
We have shown in
Section~\ref{sec:degree} that the degree of $E$ is maximal, equal to
${a+b} \choose a$, both in the generic case and in the two-step case.
This can be extended to all set steps such that $K(u)$ has at least
two (algebraically independent) indeterminate coefficients,  using the
results of~\cite{cohen}.

It would be interesting to study more cases, in particular those
involving a symmetry, which reduces the degree. Assume $\cS=-\cS$, and
$\om_{-s}=\om_s$ for all $s \in \cS$. In particular, $a=b$. 
Then, as discussed in Example~2,
the small and large roots of $K(u)$ are simply related by
$U_{a+1}=1/U_1,  \ldots, U_{2a}=1/U_a$. This implies that many products
$U_{i_1} \cdots U_{i_a}$, with $i_1<\cdots <i_a$, are actually of the
form $U_{j_1} \cdots U_{j_{a-2k}}$ for some $k>0$. The products that
reduce in that way have a minimal polynomial that strictly divides 
$$
Q(z)=\prod_{|I|=a} (1-zU_I).
$$
The non-reducing products $U_{i_1} \cdots U_{i_a}$ are the $2^a$ terms
$U_I=U_1^{\pm1} \cdots U_a^{\pm1}$. 
Thus
$$
\bar Q(z)=\prod_{\varepsilon \in \{\pm 1\}^a}(1-z U^\varepsilon),
$$
 is a polynomial in $z$ and $t$ that divides $Q(z)$, and
vanishes at $z=E$. Hence in the symmetric case, $E$ has degree at most
$2^a$. 

One could try to study systematically the cases $\cS=\llbracket -a,
a\rrbracket$ or $\cS=\{\pm1, \pm
a\}$. When $\cS=\{\pm1, \pm 2\}$, we
have seen in Example~2 that $E$ has degree 4. The Galois group $G$ of
$K(u)=u ^2-t(1+u)^2(1-u+u^2)$ over 
$\qs(t)$ can be seen to  isomorphic to the dihedral group $D_4$. More precisely,
$G=\{
\rm{id}, (1,2,3,4), (1,4,3,2), (1,3)(2,4), (1,2)(3,4),(1,4)(2,3),
(1,3), (2,4)\}$. The subgroup that leaves $U_1U_2$ invariant is the
subgroup of index 4 generated by $(1,2)(3,4)$. This explains why
$E=-U_1U_2/t$ has degree 4. 

\subsection{The \gf\ of rectangular Schur functions}\label{sec:P}
We proved in Section~\ref{sec:rectangle} that, for symmetric functions
in $n$ variables,  the \gf\ of rectangular Schur functions of height
$a$ is rational:
$$
\sum_{k\ge 0} s_{k^a} z^k = \frac{P(z)}{Q(z)},
$$
where $Q(z)$ is given by~\eqref{Q-def1} and has degree $n\choose a$, while
$P(z)$ has degree ${n\choose a}-n$. We have given two expressions
of $P(z)$ in terms of the $u_i$'s (see~(\ref{P-expr}--\ref{P-expr1})),
and proved that $P(1/u_J)$ has a  
simple expression~\eqref{PuJ}. However, we have no expansion of
$P(z)$ in 
symmetric functions, other than
$$
P(z)= \sum_{i=0}^{{n\choose a}-n}z^i \sum_{j+ k= i} (-1)^{j}e_{j}[e_a]
\ s_{k^a},
$$
which comes directly from the fact that $P(z)=Q(z) \sum_{k\ge 0}
s_{k^a} z^k$. It would be interesting to find a simpler expression for
the coefficients of $P(z)$. The term $(-1)^j$, in particular, leaves
hope for possible simplifications, which may in turn allow us to compute $P(z)$
more efficiently. Let us give the expression of $P$
for a few values of $a$ and $n$:
for $a=2$ and $ n=4$,
$$P(z)=1-e_4z^2=1-s_4z^2.$$
 For $a=2$ and $ n=5$, 
$$
P(z)= 1-{ e_4}\,{z}^{2}+{  e_{5,1}}\,{z}^{3}-{{  
    e_5}}^{2}{z}^{5}=
1-s_{{1^4}}{z}^{2}+s_{{2,1^4}}{z}^{3}-s_{{2^5}}{z}^
{5}
.$$
For $a=3$ and $n=6$,
\begin{multline*}
  P(z)=
1-s_{{2 1^4}}{z}^{2}
+ \left( s_{{2^41}}+s_{{321^4}} \right) {z}^{3}
-s_{{3^22^21^2}}{z}^{4}
- \left( s_{{3^5}}+s_{{52^5}} \right) {z}^{5}
\\
+\left( s_{{53^32^2}}+s_{{4^32^3}} +s_{{4^23^31}} \right) {z}^{6}
-2\, s_{{54^23^22}}{z}^{7}
+ \left( s_{{5^33^3}}+s_{{64^33^2}}+s_{{5^24^32}} \right){z}^{8}
\\
- \left( s_{{5^52}}+s_{{74^5}} \right) {z}^{9}
-s_{{6 ^25^24^2}}{z}^{10}
+ \left( s_{{6^454}}+s_{{765^4}} \right) {z}^{11}
-s_{{76^45}}{z}^{12}
+s_{{7^6}}{z}^{14}.
\end{multline*}
We have used the Schur basis rather than the elementary basis because
it seems, from these examples, that the coefficient of $z^i$
in $P(z)$ is either Schur-positive or Schur-negative. The conversions
to Schur functions have been made with the
package ACE 
\cite{ace}.


\subsection{The height of  random excursions} Pick an excursion
of length $n$,  uniformly at random. Its height is a random variable
$H_n$. For Dyck paths, it is known that $H_{2n}/\sqrt n$ converges in
distribution to a theta law~\cite{debruijn}, which also describes the height of the 
Brownian excursion~\cite{kennedy}. We expect this to hold for our general
excursions as well. One indication of the \emm universality, of this
law can be given in terms of trees. Via a simple bijection, 
 the height of a random Dyck path of length $2n$ becomes the
height of a random plane tree with $n$ edges. It has been proved that
the height of  other varieties of trees, like
simple trees, also follows, in the limit, a theta
law~\cite{flajolet-binary-trees}. 

The average height of plane trees
--- that is, of Dyck paths --- was derived
in~\cite{debruijn} from an expression of $E^{(k)}$ that is equivalent to our
Schur expression of this series. Is it possible, using the asymptotic
tools developed in~\cite{banderier-flajolet} for unbounded excursions,
to work out the law of the height of general excursions by starting
from our Schur expression of $E^{(k)}$?

\bigskip
\noindent
{\bf Acknowledgements.} I am indebted to Alain Lascoux and Christophe
Carré for providing efficient programs to compute plethysms.  I also
thank Alain Salinier for interesting e-discussions about the Galois
groups of trinomials, and for pointing Reference~\cite{schinzel}.

\bibliographystyle{plain}
\bibliography{biblio}

\end{document}